\documentclass[12pt]{article} 
\usepackage{amsmath, amssymb} 
\usepackage{amsthm} 
\usepackage{enumerate, url} 
\usepackage{graphicx}

\newtheoremstyle{plainsl}%
	{\topsep}
	{\topsep}
	{\slshape} 
	{}
	{\normalfont\bfseries}
	{.}
	{ }
	{}

\swapnumbers

{\theoremstyle{plainsl}
\newtheorem{theorem}{Theorem}[section]
\newtheorem{lemma}[theorem]{Lemma}
\newtheorem{corollary}[theorem]{Corollary}}
{\theoremstyle{remark}
}

\newcommand\cref[1]{Corollary~\ref{cor:#1}}
 
\newcommand\fref[1]{Figure~\ref{fig:#1}}

\renewcommand\proof{\noindent\textsl{Proof. }}
\newcommand\sqr[2]{{\vbox{\hrule height.#2pt
    \hbox{\vrule width.#2pt height#1pt \kern#1pt
        \vrule width.#2pt}\hrule height.#2pt}}}
\renewcommand\qed{%
	\ifmmode\eqno\sqr53
	\else\nolinebreak\ \hfill\sqr53\medbreak\fi}

\numberwithin{equation}{section}

\newcommand\de{\delta}
\newcommand\De{\Delta}
\newcommand\eps{\epsilon}
\newcommand\ga{\gamma}

\newcommand\sg{\sigma}

\renewcommand\th{\theta} 

\newcommand\cP{{\mathcal P}}

\newcommand\cx{{\mathbb C}}

\newcommand\ints{{\mathbb Z}}
\newcommand\re{{\mathbb R}}
\newcommand\rats{{\mathbb Q}}

\newcommand\diff{\mathbin{\mkern-1.5mu\setminus\mkern-1.5mu}}

\newcommand\seq[3]{#1_{#2},\ldots,#1_{#3}}

\newcommand\pmat[1]{\begin{pmatrix} #1 \end{pmatrix}}
\newcommand\sm[3]{\sum_{#1=#2}^{#3}}
\DeclareMathOperator{\rk}{rk}
\DeclareMathOperator{\tr}{tr}
\DeclareMathOperator{\col}{col}

\title{Periodic Graphs} 

\author{
	Chris Godsil\\
	Combinatorics \& Optimization\\
	University of Waterloo\\[2pt]
	\texttt{cgodsil@uwaterloo.ca}}

\begin{document}
\maketitle
	
\begin{abstract}
	Let $X$ be a graph on $n$ vertices with with adjacency matrix $A$ and 
	let $H(t)$ denote the matrix-valued function $\exp(iAt)$.  If $u$ and $v$ are distinct
	vertices in $X$, we say \textsl{perfect state transfer} from $u$ to $v$
	occurs if there is a time $\tau$ such that $|H(\tau)_{u,v}|=1$. If $u\in V(X)$ and
	there is a time $\sg$ such that $|H(\sg)_{u,u}|=1$, we say $X$ is
	\textsl{periodic at $u$} with period $\sg$. It not difficult to show that if 
	the ratio of distinct non-zero eigenvalues of $X$ is always rational, then $X$ is
	periodic. We show that the converse holds, from which it follows that a regular graph
	is periodic if and only if its eigenvalues are distinct.
	For a class of graphs $X$ including all vertex-transitive graphs we prove that,
	if perfect state transfer occurs at time $\tau$, then $H(\tau)$ is a 
	scalar multiple of a permutation matrix of order two with no fixed points.
	Using certain Hadamard matrices, we construct a new infinite family of graphs
	on which perfect state transfer occurs.
\end{abstract}

\section{Introduction}

Let $X$ be a graph with adjacency matrix $A$.  We define the matrix function $H(t)$
by
\[
H(t) :=\exp(itA) :=\sum_{n\ge0} i^nA^n\frac{t^n}{n!}.
\]
We note that $t$ is a real variable
\[
H(t)^* =\exp(-iAt) =H(t)^{-1}
\]
and therefore $H(t)$ is a unitary matrix.  We have $H(0)=I$ and 
\[
H(s+t)=H(s)H(t). 
\]
We say that graph is \textsl{periodic with respect to the vector $z$} if there is 
a real number $\tau$ such that $H(\tau)z$ is a scalar multiple of $z$.  (Since
$H(t)$ is unitary, this scalar will have absolute value 1.)  
We say that $X$ is \textsl{periodic relative to the vertex $u$} if it is 
periodic relative
to the standard basis vector $e_u$ or, equivalently if there is a time $\tau$ such
that $|H(\tau)_{u,u}|=1$.  We say that $X$ itself is \textsl{periodic} if there
is a time $\tau$ such that $H(\tau)$ is diagonal.

If there are basis vectors $\seq z1n$ such that $X$ is periodic with period 
$\tau_r$ relative to the vector $z_r$ for $r=1,\ldots,n$, then it follows that 
$X$ is periodic, 
with period dividing the product of the $\tau_i$'s.

We say we have \textsl{perfect state transfer} if there are distinct vertices $u$ and 
$v$ in $X$ and a time $\tau$ such that
\[
|H(\tau)_{u,v}| = 1.
\]
(If this holds then $H(\tau)_{u,v}$ is the only non-zero entry in its row 
and column.) 
Christandl, Datta, Dorlas, Ekert, Kay and Landahl \cite{christandl-2005-71} 
prove that 
we have perfect state transfer between the end-vertices of paths of
length one and two, and between vertices at maximal distance in Cartesian powers
of these graphs.  A graph exhibiting perfect state transfer between two 
vertices $u$ and $v$ models a network of quantum particles with fixed couplings, in
which the state of the particle in $u$ can be transferred to the
particle in $v$ without any information loss.

This paper is an attempt to extend some of the results in the above paper, 
and in the 
more recent work by Saxena, Severini and Shparlinski \cite{saxena-2007}.  
We prove that
a regular graph is periodic if and only if its eigenvalues are integers.  We also
show for a wide class of regular graphs (including all vertex-transitive graphs
and all distance-regular graphs) that if perfect state transfer occurs, then there
is a time $\tau$ such that $H(\tau)$ is a scalar multiple of a permutation matrix 
of order two
with zero diagonal.  (And hence the number of vertices in the graph must be even.)
Finally we present a new infinite class of antipodal distance-regular graphs where
perfect state transfer occurs.

\section{Eigenvalues}
\label{sec:evals}

We are going to derive consequences for the eigenvalues of $A$ when $X$ is periodic,
but first we show that periodicity and perfect state transfer are related.

\begin{lemma}\label{lem:}
    If perfect state transfer from $u$ to $v$ occurs at time $\tau$, then
    $X$ is periodic at $u$ and $v$ with period $2\tau$.
\end{lemma}

\proof
Assume $H=H(\tau)$ and suppose that
\[
    H(\tau)e_u = \ga e_v,
\]
where $\|\ga\|=1$. Then $H_{v,u}=\ga$ (because $H(t)$ is symmetric)
and therefore $He_v=\ga e_u$. Now
\[
    \ga^2 e_u =\ga He_v =H^2 e_u
\]
and similarly $H^2e_v=e_v$.\qed

Our main tool in this paper is the spectral decomposition of a symmetric matrix.
Suppose $A$ is a symmetric matrix with eigenvalues $\seq\th1d$, and let $E_r$
be the matrix that represents orthogonal projection on the eigenspace associated
with $\th_r$.  If $f$ a complex-valued function defined on the eigenvalues of $A$,
we have
\[
	f(A) =\sum_{r=1}^d f(\th_r)E_r.
\]
In particular,
\[
	H(t) =\sum_r \exp(it\th_r)E_r.
\]
Note that the scalars $f(\th_r)$ are the eigenvalues of $f(A)$ and that $\sum E_r=I$.
We will also make use of the fact that $f(A)$ is necessarily a polynomial in $A$
(because each projection $E_r$ is).

Assume now that $A$ is the adjacency matrix of a graph $X$ on $v$ vertices.
If $x$ is a vector in $\re^v$, we say that an eigenvalue $\th_r$ of $X$ is
in the \textsl{eigenvalue support} of $x$ if $E_rx\ne0$.  
The case of interest to
us will be when $x$ is the characteristic vector of a subset $S$ of $V(X)$.
The eigenvalue support of a non-zero vector cannot be empty and, if $x\ne0$,
we define the \textsl{dual degree} of $x$ to be one less than the size of its
eigenvalue support.

\begin{theorem}
\label{thm:ratio}
	Let $X$ be a graph and let $u$ be a vertex in $X$ at which $X$ is periodic.
	If $\th_k$, $\th_\ell$, $\th_r$, $\th_s$ are eigenvalues in the support
	of $e_u$ and $\th_r\ne\th_s$,
	\[
		\frac{\th_k-\th_\ell}{\th_r-\th_s}\in\rats.
	\]
\end{theorem}

\proof
Suppose the period of $X$ at $1$ is $\tau$.  Then
\[
	H(\tau)_{1,1} =\ga
\]
where $|\ga|=1$.  Since $H(\tau)$ is unitary, this implies that we have
\[
	H(\tau) =\pmat{\ga&0 \\0&H_1}
\]
where $H_1$ is unitary.  Hence
\[
	\ga e_1 = H(\tau)e_1 =\sum_r \exp(i\tau\th_r) E_r e_1
\]
The non-zero vectors $E_re_1$ are pairwise orthogonal and so linearly 
independent, whence we see that
\[
	\exp(i\tau\th_r) =\ga
\]
for each eigenvalue $\th_r$ in the support of $e_1$.

Consequently of $\th_r$ and $\th_s$ are in the eigenvalue support of $e_1$, then
\[
	\exp(i\tau(\th_r-\th_s)) =1
\]
and so $\tau(\th_r-\th_s)$ is an integral multiple of $2\pi$.  This proves the theorem.\qed

\begin{corollary}
	If $X$ is periodic at $u$ and there are two integer eigenvalues in the support 
	of $e_u$, then all eigenvalues in its support are integers.
\end{corollary}

\proof
Suppose $k$ and $\ell$ are distinct integer eigenvalues in the support of $e_u$.  By
the theorem
\[
	\frac{\th_r-\ell}{k-\ell}\in\rats
\]
and therefore $\th_r$ is rational. Since $\th_r$ is an algebraic integer, it must
be an integer.\qed

The condition in the previous theorem is known as the \textsl{ratio condition}.  
Christandl et al \cite{christandl-2005-71}
stated that if perfect state transfer takes place on a path, then its eigenvalues must
satisfy the ratio condition.  Similarly Saxena et al \cite{saxena-2007} proved that 
if a circulant graph is periodic, its eigenvalues must satisfy the ratio condition.

We can derive these results from ours, taking the second first.  A circulant
graph is vertex transitive and so if $E_r e_u=0$ for some vertex $u$, then
for all vertices $u$, we have $E_r e_u=0$.  But this implies that $E_r=0$,
and so we conclude that the support of $e_u$ is the set of all eigenvalues of $X$.
Hence the eigenvalues of a periodic circulant satisfy the ratio condition.

Now paths.  From the spectral decomposition we have
\[
	(tI-A)^{-1}_{u,u} =\sum_r (t-\th_r)^{-1}(E_r)_{u,u}
\]
while standard matrix theory yields that
\[
	(tI-A)^{-1}_{u,u} = \frac{\phi(X\diff u,t)}{\phi(X,t)}.
\]
Since
\[
	(E_r)_{u,u} = e_u^T E_r e_u =e_u^T E_r^2e_u =e_u^TE_r^TE_re_u,
\]
we see that $(E_r)_{u,u}$ is not zero if and only if $E_re_u$ is not zero, that
is, if and only if $\th_r$ is in the support of $e_u$.  
Therefore the dual degree of $e_u$ is
equal to the number of poles of the rational function $\phi(X\diff u,t)/\phi(X,t)$,
less 1.  The characteristic polynomial of the path $P_n$ on $n$ vertices
satisfies the recurrence
\[
	\phi(P_{n+1},t) =\phi(P_n,t) - \phi(P_{n-1},t),\qquad (n\ge1);
\]
since $\phi(P_0,t)=1$ and $\phi(P_1,t)=t$, it follows that $\phi(P_n,t)$ and
$\phi(P_{n-1},t)$ are coprime.  So our rational function has $n$ distinct zeros,
and again we find that the eigenvalue support of a vertex of degree one in a path
is the set of all eigenvalues of the path.

\section{Periodicity}

It follows from the spectral decomposition
\[
	H(\tau) =\sum_r \exp(i\tau\th_r) E_r
\]
that if the eigenvalues of $X$ are integers, then 
\[
	H(2\pi) =\sum_r E_r =I
\]
and therefore $X$ is periodic with period $2\pi$.  This shows that the path $P_2$
is periodic.  However, as Christandl et al noted, $P_3$ is periodic and its
eigenvalues are
\[
	-\sqrt2,0,\sqrt2.
\]
Thus integrality is sufficient, but not necessary.  There is a more general
condition which is necessary and sufficient.  

\begin{theorem}\label{thm:period}
	Suppose that the eigenvalues of $X$ are $\seq\th1d$.  Then $X$ is periodic 
	if and only if the ratio of any two non-zero eigenvalues is rational.
\end{theorem}

\proof
We first show that the given condition is sufficient.  If it holds and $\th_s\ne0$, 
there are rational numbers $\seq{q}1r$ such that $\th_r= q_r\th_s$ and 
therefore there is an integer $N$ and integers $\seq{k}1r$ such that $N\th_r =k_r\th_s$.  
Hence if we take
\[
	\tau = \frac{2N\pi}{\th_s}
\]
we have
\[	
	\exp(i\tau\th_r) =\exp(2ik_r\pi) =1
\]
for all $r$, and consequently $X$ is periodic.

Conversely, suppose $X$ is periodic with period $\tau$.
Then $H(\tau)$ is diagonal and commutes with $A$. If we view the diagonal 
of $H(\tau)$ as a function on $V(X)$, it must be constant on the connected
components of $X$. 

If $X$ is connected it follows that there is a complex 
number $\ga$ such that $|\ga|=1$ and $H(\tau)=\ga I$.  
Since the sum of the eigenvalues of $X$ is zero,
the product of the eigenvalues of $H(\tau)$ is one, and thus if
$v=|V(X)|$, then
\[
	1 = \det(H(\tau)) =\ga^v.
\]
Therefore $H(v\tau)=I$, and so all eigenvalues of $H(v\tau)$ are equal to one.
If $X$ is not connected, we can apply this argument to deduce that each
diagonal entry of $H(\tau)$ is an $m$-th root of unity, for some integer $m$,
and therefore if $w$ is the least common multiple of the sizes of the
connected components of $X$, then $H(w\tau)=I.$

Hence
\[
	\exp(iw\tau\th_r) =1
\]
for all $r$, and there are integers $\seq\ell1r$ such that
\[
	w\tau\th_r =2\ell_r\pi.
\]
This implies that, for all $r$ and $s$,
\[
	\frac{\th_r}{\th_s} = \frac{\ell_r}{\ell_s} \in\rats.\qed
\]

\begin{lemma}\label{lem:}
    If the ratios of the non-zero eigenvalues of $X$ are rational, then the
    square of an eigenvalue is an integer.
\end{lemma}

\proof
If the condition of the lemma holds and $\th_s\ne0$, there are rationals 
$\seq{q}1r$ such that $\th_r =q_r\th_s$. 
Since the sum of the squares of the eigenvalues of $A$
is equal to twice the number of edges, we see that $\th_s^2$ is rational.
Since $\th_s$ is an algebraic integer, $\th_s^2$ is thus a rational algebraic
integer, and accordingly it is actually an integer.\qed

\begin{corollary}\label{cor:period}
    A graph is $X$ is periodic if and only if either:
    \begin{enumerate}[(a)]
        \item 
        The eigenvalues of $X$ are integers, or
        \item
        The eigenvalues of $X$ are rational multiples
        of $\sqrt\De$, for some square-free integer $\De$.
    \end{enumerate}
    If the second alternative holds, $X$ is bipartite.
\end{corollary}

\proof
The stated conditions are sufficient, and so we show they are necessary.

Suppose $k$ is a non-zero integer eigenvalue of $X$.  Then $\th_r/k\in\rats$
and so  $\th_r\in\rats$.  Since $\th_r$ is an algebraic integer, it is therefore
an integer. If no non-zero eigenvalue of $X$ is an integer, then 
by the theorem each eigenvalue
is a rational multiple of the spectral radius $\th_1$, say $\th_i = m_i\th_1$.

Also by the theorem $\th_1^2\in\ints$, whence it follows that $t^2-\th_1^2$ is
the minimal polynomial of $\th_1$ over the rationals, and therefore $-\th_1$
is an eigenvalue of $X$. By the Perron-Frobenius theory, this implies that $X$
is bipartite.\qed

\section{Coherent Algebras}

A \textsl{coherent algebra} is a real or complex algebra of matrices that 
contains the all-ones matrix $J$ and is closed under Schur multiplication,
transpose and complex conjugation.  (The \textsl{Schur product} $A\circ B$ 
of two matrices $A$ and $B$ with same order is defined by
\[
(A\circ B)_{i,j} = A_{i,j}B_{i,j}.
\]
It has been referred to as the ``bad student's product''.) 
A coherent algebra always has a basis of $01$-matrices and this is unique, given
that its elements are $01$-matrices.  This set of matrices determines a set of
directed graphs and the combinatorial structure they form is known as a
\textsl{coherent configuration}.  When we say that a graph $X$ is a 
``graph in a coherent algebra'', we mean that $A(X)$ is a sum of distinct 
elements of the $01$-basis.

A coherent algebra is \textsl{homogeneous} if the identity matrix is an element
of its $01$-basis. If $M$ belongs to a homogeneous coherent algebra, then
\[
M\circ I = \mu I
\]
for some scalar $\mu$. Hence the diagonal of any matrix in the algebra is
constant. If $A$ is a $01$-matrix in the algebra, the diagonal entries of $AA^T$
are the row sums of $A$. Therefore all row sums and all column sums of any matrix
in the $01$-basis are the same, and therefore this holds for each matrix in the
algebra. In particular we can view the non-identity matrices as adjacency
matrices of regular directed graphs. Any directed graph in a homogeneous coherent
algebra must be regular.

We consider two classes of examples.  First, if $\cP$ is a set of permutation
matrices of order $n\times n$, then the commutant of $\cP$ in 
$\textrm{Mat}_{n\times n}(\cx)$
is Schur-closed.  
Therefore it is a coherent algebra, and this algebra is homogeneous
if and only the permutation group generated by $\cP$ is transitive.
Thus any graph whose automorphism group acts transitively on its vertices
belongs to a coherent algebra

For second class of examples, if $X$ is a graph we define the \textsl{$r$-th
distance graph} $X_r$ of $X$ to be the graph with $V(X_r)=V(X)$, where vertices
$u$ and $v$ are adjacent in $X_r$ if and only if they are at distance $r$ in $X$.
So if $X$ has diameter $d$, we have distance graphs $\seq X1d$ with corresponding
adjacency matrices $\seq A1d$. If we use $A_0$ to denote $I$, then the graph $X$
is \textsl{distance regular} if the matrices $\seq A0d$ are the $01$-basis of a
homogeneous coherent algebra. (It must be admitted that this is not the standard
definition.) We will refer to the matrices $A_i$ as \textsl{distance matrices}.

A commutative coherent algebra is the same thing as a Bose-Mesner algebra of an
association scheme. We will not go into this here, but we do note that the
coherent algebra belonging to a distance-regular graph is commutative.

\begin{theorem}\label{thm:}
	If $X$ is a graph in a coherent algebra with vertices $u$ and $v$, and 
	perfect
	state transfer from $u$ to $v$ occurs at time $\tau$, then $H(\tau)$ is a 
	scalar multiple of a permutation
	matrix with order two and no fixed points that lies in the centre of the
	automorphism group of $X$.
\end{theorem}

\proof
First, if $A=A(X)$ where $X$ is a graph in a homogeneous coherent algebra then,
because it is a polynomial in $A$, the matrix $H(t)$ lines in the algebra for all
$t$.
Hence if
\[
|H(\tau)_{u,v}| = 1
\]
it follows that $H(\tau) =\xi P$ for some complex number $\xi$ such that $|\xi|=1$
and some permutation matrix $P$.  Since $A$ is symmetric, so is $H(t)$ for any $t$,
and therefore $P$ is symmetric.  So
\[
P^2 =PP^T =I
\]
and $P$ has order two.  Since $P$ has constant diagonal, its diagonal is zero
and it has no fixed points.  As $P$ is a polynomial in $A$, it commutes with
any automorphism of $X$ and hence is central.\qed

\begin{corollary}\label{cor:evncoh}
	If $X$ is a graph in a coherent algebra with vertices $u$ and $v$ and perfect
	state transfer from $u$ to $v$ occurs at some time, then the number of vertices 
	of $X$ is even.
\end{corollary}

\proof
Since $P^2=I$ and the diagonal of $P$ is zero, the number of columns of $P$ is even.\qed

Saxena et al \cite{saxena-2007} proved this corollary for circulant graphs.  
A homogeneous coherent
algebra is \textsl{imprimitive} if there is a non-identity matrix in its $01$-basis
whose graph is no connected, otherwise it is \textsl{primitive}.  If the algebra
is the commutant of a transitive permutation group, it is imprimitive if and
only the group is imprimitive as a permutation group.  The above corollary implies that 
if perfect state transfer takes place on a graph from a homogeneous coherent algebra,
the algebra is imprimitive.

Note that our corollary holds for any vertex-transitive graph, and for any distance-regular
graph.

\section{Walk-Regular Graphs}

A graph is \textsl{walk regular} if its vertex-deleted subgraphs $X\diff u$
are all cospectral.
We have the following result which implies that any graph in a coherent
configuration is walk regular. It follows immediately from Theorem~4.1 
in \cite{MR568777}.

\begin{lemma}\label{lem:wr-diag}
	A graph $X$ with adjacency matrix $A$ is \textsl{walk regular} 
	if and only if the diagonal
	entries of $A^k$ are constant for all non-negative integers $k$.
\end{lemma}

Any connected regular graph with at most four distinct eigenvalues is walk
regular (see Van Dam \cite{MR1344560}). Note that a walk-regular graph is
necessarily regular.
The graph in \fref{wr12a} is walk regular but not vertex transitive.  
This graph does not lie in a homogeneous 
coherent algebra---the row sums of the Schur product
\[
A\circ (A^2-4I)\circ (A^2-4I-J)
\]
are not all the same. 

\begin{figure}[htbp]
	\centering
		\includegraphics[scale=0.5]{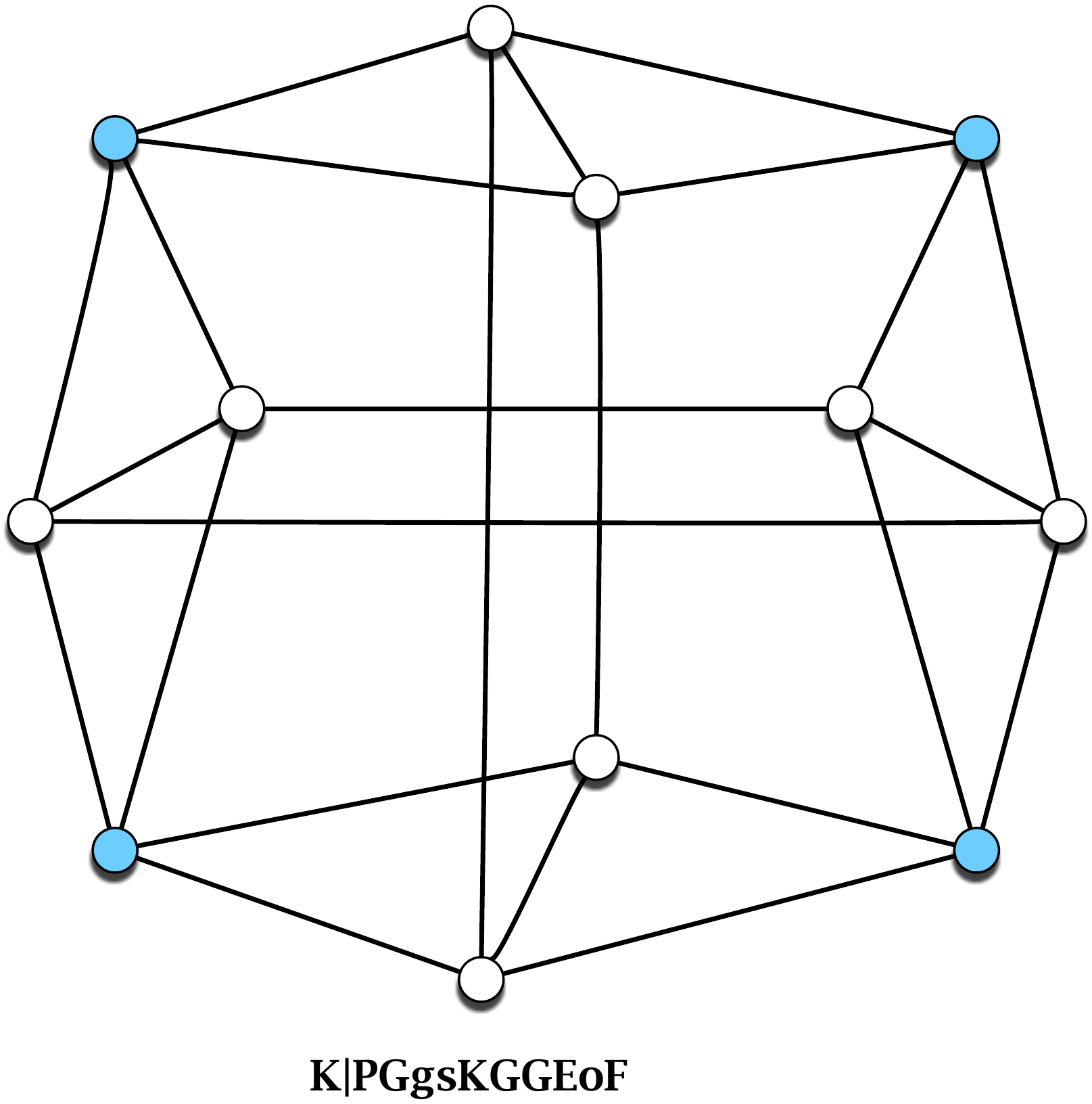}
	\caption{A Walk-Regular Graph}
	\label{fig:wr12a}
\end{figure} 

Suppose perfect state transfer takes place from $u$ to $v$ at time $\tau$.  
Then there is a complex number $\ga$ where $|\ga|=1$ such that
\[
	H(2\tau)_{u,v} =H(2\tau)_{v,u} =\ga.
\]
As $|\ga|=1$, we see that $H_{u,u}=0$ and, as $X$ is walk-regular, all diagonal
entries of $H$ are zero and in particular, $\tr(H(\tau))=0$. Since
\[
	H(2\tau)_{u,u} =\ga^2
\]
we also see that that $H(2\tau)=\ga^2I$. We conclude that the eigenvalues of
$H(\tau)$ are all $\pm\ga$ and, since $\tr(H(\tau))=0$, both $\ga$ and $-\ga$
have multiplicity $|V(X)|/2$. (Accordingly $|V(H)|$ must be even.)

If $m_\th$ denotes the multiplicity of the eigenvalue $\th$,
then 
\[
\tr(H(t)) =\sum_\th m_\th \exp(it\th)
\] 
and thus we may define the Laurent polynomial $\mu(z)$ by
\[
\mu(z) := \sum_\th m_\th z^\th.
\] 
Here $z=\exp(it)$ and we have made use of the fact that $X$ is periodic and hence
its eigenvalues are integers.  We call $\mu(z)$ the \textsl{multiplicity enumerator} of $X$.

\begin{lemma}\label{lem:}
	Let $X$ be a walk-regular graph with integer eigenvalues.  If perfect state
	transfer occurs on $X$, then $\mu(z)$ has an eigenvalue on the unit circle 
	of the complex plane.
\end{lemma}

\proof
The trace of $H(t)$ is zero if and only if $\mu(\exp(it))=0$.\qed

The characteristic polynomial of the graph in \fref{wr12a} is
\[
(t-4)(t-2)^3t^3(t+2)^5
\]
and its multiplicity enumerator is
\[
\mu(z) =z^{-2}(z^6+3z^4+3z^2+5).
\]
This polynomial has no roots on the unit circle, and we conclude that perfect
state transfer does not occur.

We see no reason to believe that, if perfect state transfer occurs on a
walk-regular graph at time $\tau$, then $H(\tau)$ must be a multiple of 
a permutation matrix. (But we do not have an example.)

For more information on walk-regular graphs, see \cite{MR568777}. The
computations in this section were carried out in sage \cite{sage}. Van Dam
\cite{MR1344560} studies regular graphs with four eigenvalues, as we noted these
provide examples of walk-regular graphs.

\section{Hypercubes}

The $n$-cube was one of the first cases where perfect state transfer was shown to
occur (see \cite{christandl-2005-71}). The $n$-cube is an example of a
distance-regular graph, and in this section we establish some of the consequences
of our work for distance-regular graphs.

Suppose $X$ is a distance-regular graph with diameter $d$ and distance matrices
$\seq A0d$. It has long been known that if $X$ is imprimitive then either $X_2$
is not connected or $X_d$ is not connected (see \cite[Chapter~22]{nbagt} or
\cite[Section~4.2]{bcn}). For the $n$-cube, neither $X_2$ nor $X_n$ are
connected. If $X_2$ is not connected then $X$ is bipartite and $X_m$ is not
connected if $m$ is even. Here we will be concerned with the case where $X_d$ is
not connected. In this case it follows from the work just cited that the
components of $X_d$ must be complete graphs of the same size, and that two
vertices of $X$ lie in the same component of $X_d$ if and only if they are at
distance $d$ in $X$. The components are called the \textsl{antipodal classes} of
$X$; note these may have size greater than two. If $X_d$ is not connected, we say
that $X$ is an \textsl{antipodal} distance regular graph.

As already noted the $n$-cubes provide natural examples of distance-regular
graphs. The $n$-cube has diameter $n$, and its vertices can be partitioned into
pairs such that two vertices in the same pair are at distance $n$, while vertices
in different pairs are at distance less than $n$, hence it is antipodal.
 
Our next lemma is more or less a specialization of \cref{evncoh} to the case of
distance-regular graphs.

\begin{lemma}\label{lem:pfsdrg}
	Suppose $X$ is a distance-regular graph with diameter $d$.  If perfect state 
	transfer from $u$ to $v$ occurs at time $\tau$, then $X$ is antipodal with 
	antipodal classes of size two and $H(\tau)=X_d$.\qed
\end{lemma}

We outline a proof that perfect state transfer occurs on the $n$-cube.  In the following
section we will use much the same argument to provide a second infinite class of examples.

The $n$-cube gives rise to an association scheme, that is, a commutative coherent
configuration as follows. Let $X$ be the graph of the $n$-cube, and view its
vertices as binary vectors of length $n$. Let $X_i$ be the graph with the same
vertex set as $X$, where two vertices are adjacent if and only if they differ in
exactly $i$ positions. Thus $X_1=X$ and two vertices are adjacent in $X_i$ if and
only they are at distance $i$ in $X$. Let $A_i$ be the adjacency matrix of $X_i$
and set $A_0=I$. Then $\seq A0n$ are symmetric $01$-matrices and
\[
\sm r0n A_r =J.
\]
It is known that these matrices commute and that their span is closed under
matrix multiplication. Hence they form a basis for a commutative coherent algebra
of dimension $d+1$. Since this algebra is commutative and semisimple, it has a
second basis
\[
\seq E0n
\]
of matrix idempotents.  These are real symmetric matrices such that 
$E_iE_j =\de_{i,j}$ and
\[
 \sm i0n E_i =I.
\]
These idempotents represent orthogonal projection onto the eigenspaces of 
the algebra and hence there are scalars $p_i(j)$ such that 
$A_iE_j =p_i(j)E_j$.  Consequently
\[
A_i =A_i I =A_i\sm r0n E_r =\sm r0n p_i(r)E_r.
\]
The scalars $p_i(j)$ are known as the \textsl{eigenvalues} of the scheme.

On the other hand, the matrices $\seq A0n$ are a basis for the algebra and thus
there are scalars $q_i(j)$ such that
\[
E_j =2^{-n}\sm r0n q_j(r)A_r.
\]
The scalars $q_i(j)$ are the \textsl{dual eigenvalues} of the scheme.

If $f$ is any function defined on the eigenvalues of the scheme then
\[
f(A_1) =\sm r0n f(p_1(r))E_r.
\]
For the $n$-cube, $p_1(r)=n-2r$ and accordingly
\[
H(t) = \sm r0n \exp(i(n-2r)t)\,E_r =\exp(int) \sm r0n \exp(-2irt)E_r.
\]
If we set $t=\pi$, we then have
\[
H(\pi) =(-1)^n \sm r0n E_r = (-1)^nI.
\]
Thus the $n$-cube is periodic with period $\pi$.  Now try $n=\pi/2$.  Then
\[
H(\pi/2) =i^n \sm r0n (-1)^r E_r
\]
From, for example, \cite[Section~21.3]{MacSlo} we know that
\[
\sm r0n (-1)^r E_r =A_d
\]
and we conclude that perfect state transfer occurs at time $\pi/2$.

In \cite{jafarizadeh-2007}, Jafarizadeh and Sufiani give examples of perfect state transfer 
on distance-regular graphs that arise when a weighted adjacency matrix is used
in place of the ordinary adjacency matrix.

\section{Hadamard Matrices}

We now show how to use a special class of Hadamard matrices to construct
distance-regular graphs of diameter three with perfect state transfer.

The Hadamard matrices $H$ we want must be as follows:
\begin{enumerate}[(a)]
	\item 
	Symmetric: $H=H^T$.
	\item
	Regular: all row and column sums are equal.
	\item
	Constant diagonal: all entries on the diagonal are equal.
\end{enumerate}
Suppose $H$ is $n\times n$ and satisfies this constellation of conditions.
If its rows sum to $c$, then $HJ=cJ$ and $H^2J=c^2J$. Since
\[
    H^2 I =HH^T =nI,
\]
whence $c=\pm\sqrt{n}$.  It follows that there are two possibilities:
either the diagonal entries and the rows sums of $H$ have the same
sign, or they do not.  To handle this we assign a parameter $\eps=\eps(H)$
to $H$, by defining $\eps$ to be the sign of the product of the row sum
of $H$ with a diagonal entry.  Whatever the value of $\eps(H)$, the order of
$H$ must be a square.  For a recent useful paper on these Hadamard matrices,
see Haemers \cite{wh-maxen}.

We remark that if the row sums of a Hadamard matrix are all equal then its columns
sums must all be equal as well, and again the order of the matrix must be a square.
Similarly, a symmetric Hadamard matrix with constant diagonal must have square order.
If $H$ and $K$ are regular symmetric Hadamard matrices with constant diagonal,
then so is their Kronecker product $H\otimes K$ and 
\[
\eps(H\otimes K) =\eps(H)\eps(H).
\]
Finally if $H$ is a regular Hadamard matrix of order $n\times n$ with constant diagonal, 
and $P$ is the permutation matrix on $\re^n\otimes\re^n$ defined by
\[
P(u\otimes v) =v\otimes u,
\]
then $P(H\otimes H^T)$ is a symmetric regular Hadamard matrix with constant diagonal.

Suppose $H$ is a regular symmetric $n\times n$ Hadamard matrix with constant diagonal.
Replacing $H$ by $-H$ is necessary, we may assume that the diagonal entries
of $H$ are equal to $1$, and then the matrix
\[
\frac12(J+H)-I
\]
is the adjacency matrix of a regular graph (in fact a strongly regular graph).
Unfortunately the graphs that result are not the ones we need.
Our aim rather is to construct an antipodal distance-regular graph of diameter three
on $2n$ vertices from $H$.  (These graphs are member of a class known as 
\textsl{regular two-graphs}.  An introductory treatment is offered 
in \cite[Chapter~11]{MR1829620})

We construct our graphs as follows.  Let $H$ be a symmetric Hadamard matrix
of order $n\times n$ with constant row sum and with all diagonal entries
equal to 1.  We construct a graph $X(H)$ with vertex set
\[
\{1,\ldots,n\} \times\{0,1\},
\]
where $(i,a)$ is adjacent to $(j,b)$ if and only if $i\ne j$ and
\[
H_{i,j} = (-1)^{a+b}  
\]
Clearly this gives us a regular graph $X(H)$ on $2n$ vertices with valency $n-1$.
It is not hard to check that its diameter is three, and two vertices are
at distance three if and only if they are of the form 
\[
(i,0),\ (i,1)
\]
for some $i$.  Suppose $A$ is the adjacency matrix of $X(H)$.  The space
$S$ of vectors that are constant on the antipodal classes of $X(H)$ has dimension $n$
and the matrix representing the action of $A$ on $S$ is the $n\times n$ matrix $J_I$.  
Its eigenvalues are $n-1$ and $-1$ with respective multiplicities 1 and $n-1$.  
The orthogonal complement of $S$ is also $A$-invariant and the matrix representing 
the action of $A$ is $H$.  

If we write $H$ as $I+C-D$ where $C$ and $D$ are 01-matrices, then
\[
A(X(H)) =C\otimes\pmat{1&0\\0&1} +D\otimes\pmat{0&1\\1&0}
\]
and the eigenvalues of $A$ are the eigenvalues of the matrices
\[
C+D =J-I,\quad C-D =H-I.
\]
Therefore the eigenvalues of $A$ are $n-1$, $-1$ and the eigenvalues
of $H-I$. Since $H^2=nI$, we see that the eigenvalues of $A$
are
\[
n-1,\ -1,\ \sqrt{n}-1,\ -\sqrt{n}-1.
\]  
If $a$ is the multiplicity of $\sqrt{n}$ as an eigenvalue of $F$, then
\[
n =\tr(A) =a\sqrt{n} -(n-a)\sqrt{n}
\]
whence
\[
a = \frac12(n+\sqrt{n}).
\]
it follows that
\begin{align*}
\mu(z) &= z^{n-1} +(n-1)z^{-1} +\frac12(n+\sqrt{n}) z^{\sqrt{n}-1}
									+\frac12(n-\sqrt{n}) z^{-\sqrt{n}-1}\\
       &=z^{-\sqrt{n-1}} (z^{n+\sqrt{n}} +(n-1)z^{\sqrt{n}}
				+\frac12(n-\sqrt{n})z^{2\sqrt{n}}+\frac12(n-\sqrt{n})).
\end{align*} 
If $z^{\sqrt{n}}=-1$, then $z^{2\sqrt{n}}=1$ and, since $\sqrt{n}$ is even,
\[
    z^{n-\sqrt{n}} =(z^{\sqrt{n}})^{\sqrt{n}} z^{-\sqrt{n}} =-1.
\]
Therefore if $z^{\sqrt{n}}=-1$, then $\mu(z)=0$ and if $z=e^{it}$ then
\[
    H(t) = -E_{n-1} -E_{(-1)} +E_{(-1+\sqrt{n})} +E_{(-1-\sqrt{n})}.
\]

Let $F$ be the direct sum of $n$ copies of the matrix
\[
    \frac12\pmat{1&1\\ 1&1}.
\]
Then $\rk(F)=n$ and $\col(F)$ is the space of functions constant on the antipodal classes
of $X(H)$.  Since $F^2=F$, we find that $E_{(-1)}=F-E_{n-1}$ and consequently
\[
    H(t) =I -2F.
\]
So $H(t)$ is a permutation matrix and we have perfect state transfer.

Infinitely many examples of regular symmetric Hadamard matrices with constant diagonal
are known, in particular they exist whenever the order $n$ is a power of four.
Haemers \cite{wh-maxen} provides a summary of the state of our knowledge.

\section{Questions and Comments}

The most obvious problem remaining is to characterize the pairs of vertices
where perfect state transfer takes place. 

Angeles-Canul et al \cite{acnoprt} et al give examples of graphs where perfect state 
transfer occurs but the underlying graph is not periodic; this answered a question 
we raised in an earlier version of this paper.

Facer, Twamley and Cresser \cite{facer-2007} construct certain Cayley 
graphs for $\ints_2^d$ on which perfect state transfer; the resulting 
family of graphs include the $n$-cubes. In \cite{bgs}, it is shown 
that if $C$ is the connection set
for a Cayley graph for $\ints_2^d$ and $\sg$ is the sum of the elements of $C$
and $\sg\ne0$, then we have perfect state transfer from $v$ to $v+\sg$ for 
each element $v$ of $\ints_2^d$. 
The $d$-cube has $2^d$ edges. It would be interesting to have a good estimate 
for the minimum number of edges in a graph where perfect state transfer between
two vertices at distance $d$ occurs.

I would like to acknowledge a number of very useful discussions with 
Simone Severini  on the topics of this paper.

\end{document}